# PROFILE-DRIVEN AUTOMATED MIXED PRECISION

Ralph Nathan[1], Helia Naeimi[2], Daniel J. Sorin[1], Xiaobai Sun[3]

[1]Department of ECE
Duke University

[2]Intel Labs
Intel Corporation

[3]Department of Computer Science
Duke University

**ABSTRACT**

*We present a scheme to automatically set the precision of floating point variables in an application. We design a framework that profiles applications to measure undesirable numerical behavior at the floating point operation level. We use this framework to perform mixed precision analysis to heuristically set the precision of all variables in an application based on their numerical profiles. We experimentally evaluate the mixed precision analysis to show that it can generate a range of results with different accuracy and performance characteristics.*

## 1. INTRODUCTION

Modern computer architectures allow for trade-offs to be made between accuracy and performance by supporting multiple precisions for floating point values and arithmetic operations. Higher precision improves accuracy but is often much slower than lower precision. Thus, a programmer of numerical software seeks to use only as much precision as needed to achieve the desired accuracy.

Although it may be convenient for a programmer to simply declare all operations in a program to be at a given precision (e.g., declaring all variables to be of the lower-precision type float or all variables to be of the higher-precision type double), a sophisticated programmer may wish to employ "mixed precision," in which some operations are performed at higher precision (because they are critical to accuracy) and the rest are performed at lower precision (because they are *not* as critical to accuracy). The challenge with mixed precision is determining which operations should be performed at which precision. In this work, we consider mixed precision only for *static* operations. For example, the addition performed by a given add instruction is always at the same precision during a given execution; we do not change the precision during execution.

There are two existing approaches to mixed precision programming. The first approach is *programmer-directed*; the programmer manually chooses the precision for each operation. This programmer-directed approach is feasible when the program has a relatively simple structure and its numerical behavior is well understood by the programmer. Because these conditions are rarely true, programmers are often tempted by the simple, but inefficient, solution of declaring all operations to be at higher precision.

The second approach to mixed precision is to perform an *automated search* of the mixed-precision space [1][2][3]. Automation is a key feature and it is one that we build on in our work. Prior automated schemes, which we describe in more detail in Section 7, run different mixed-precision versions of the program; that is, each version has a different assignment of static operations to precisions in a search for the best outcome. The major drawback to this search approach is the enormous scope and cost of the search. With two precision levels, a program with $N$ floating point operations has $2^N$ possible mixed precision configurations. For most programs, sampling even a small fraction of this space is extremely time-consuming.

To address the limitations of prior automated approaches to mixed precision, we propose a *profile-driven* tool, called AMP (for Automated Mixed Precision), to automatically set the precision of all floating point operations in an application. The key to AMP's ability to search the huge space of mixed precision versions is its use of profiling to narrow down the search. The input to AMP is a single-precision application in which all operations are at the minimum precision, and AMP's output is a mixed-precision application in which precisions have been chosen to improve accuracy. There are two steps to this automatic process. First, AMP profiles the application to uncover, for each floating point operation, numerical faults, i.e., any numerically ill behaviors that could degrade accuracy. Second, AMP uses this profile to identify floating point operations with numerical faults and increase their precision.

In developing AMP, we make the following contributions:

- We categorize the numerically pathological conditions that lead to significant loss of numerical accuracy.
- We show how to monitor and profile these conditions during execution.
- We show to use the profile to (a) narrow down the search in the huge space of possible mixed precision configurations and (b) use heuristics to set operation precisions accordingly.
- We experimentally show that AMP can automatically and efficiently create mixed-precision applications.

AMP, like all prior work in mixed precision, has two fundamental limitations. First, the mixed precision setting is static, with its effectiveness depending on the dataset used for profiling. Sometimes, numerical programmers have model problems and datasets for testing their algorithms and programs. Second, although we can monitor and locate numerical faults, it is infeasible to trace and quantify error propagation through every computational sequence of operations, let alone their response to various precision settings. Nevertheless, despite the intractability of finding optimal mixed precision settings, automatically finding good mixed precision settings that work well for a wide range of input datasets is both desirable and, as we show here, possible.

Table 1. Profiling: Behaviors and Thresholds

| Numerical Faults | Threshold for Promoting a Static Floating Point Operation |
| --- | --- |
| Large round-off (Sec. 3.1.1) | If more than $t_1$% of instances have more than $t_2$% Error Ratio |
| Large difference in addend exponents (Sec. 3.1.2) | If more than $t_3$% of instances have exponents that differ by more than $t_4$ |
| Severe cancellation (Sec. 3.1.3) | If *any* instance has more than $t_5$ of its most significant bits cancelled |
| Near overflow/underflow (Sec. 3.1.4) | If *any* instance has a result with magnitude greater than $t_6$ or less than $t_7$ |

## 2. MOTIVATION

Automatic detection of numerical faults and automatically setting precision settings are desirable for a wide range of algorithms, especially iterative algorithms in which numerical solutions are corrected at different accuracy levels and gain more significant bits as the iteration process converges. Consider the following two specific reasons for using lower precision for some operations. First, not all operations and variable values in an application are equally critical to numerical accuracy. Often, the accuracy-critical spots are sparse in the course of execution as well as among variables. Second, modern architectures often have much better performance on an application that uses single-precision than on the same application re-written to use double-precision. Intel datasheets reveal a 2X performance difference for CPUs [4] and greater than 2X for the latest GPUs [5]. Although we focus on CPUs in this paper, we note that similar analyses can also be performed for GPUs.

To understand the performance disparity between single and double precision, we consider scalar and vectorized applications separately. We further assume a typical modern processor in which floating point arithmetic is performed on Single Instruction Multiple Data (SIMD) hardware. Consider an *n*-bit wide SIMD unit; at any given time, it can perform either *n*/64 64-bit operations or twice as many (i.e., *n*/32) 32-bit operations.

Vectorized software is tailor-made for SIMD hardware and can use the full bandwidth of the SIMD hardware. Thus the performance of vectorized 32-bit software will be significantly greater than the performance of vectorized 64-bit software. This trade-off inspires programmers to use mixed-precision, in which each operation is performed at the precision necessary for the desired accuracy (but no greater precision).

However, if the code is not vectorized, then single and double precision are expected to have similar performance characteristics for compute-bound applications. This lack of vectorization might be due to limited data level parallelism in the algorithm or the inability of the programmer to optimize her code to maximize the utilization of the hardware. The reason why scalar compute-bound applications have roughly the same performance for single and double precision is that the *latency* (but not throughput) of most single precision floating point operations is the same as that of double precision operations.

## 3. PROFILING

AMP takes a single precision application as input and profiles it. Based on the profile, AMP promotes operations (i.e., increases their precision from single to double) with bad numerical behaviors.[1] AMP profiles applications at the floating point operation level; for each floating point operation, AMP detects and quantifies numerically ill behaviors. The numerical profile consists of the Dynamic Data Flow Graph (DDFG) annotated with the numerical behaviors discussed next.

### 3.1 Numerical Behaviors Profiled

Although some prior work [6][7][8][9][10], described more in Section 7, has monitored numerical behaviors, there did not exist an established and standard set of numerical behaviors to profile. Partly inspired by prior work, but also inspired by a desire for a simple hardware implementation, we propose a small and simple set of *numerical faults* that indicate or lead to accuracy degradation.

AMP profiles the rounding error that occurs after each floating point operation, as well as a set of numerical behaviors proposed across several prior papers. We summarize these behaviors, as well as the thresholds used to determine when those behaviors are harmful, in Table 1.

Although we prototyped this profiler in software, we designed it so that it could be practically and economically implemented in hardware. Specifically, we make approximations in our profiles to facilitate hardware implementation.

---

[1] We could start with a double precision application, downgrade all operations to single precision, and then apply AMP to this single precision application.

### 3.1.1 Large Round-off

In any floating point operation, the rounding error captures all information lost due to limited precision. For any operation, ∘, with inputs $a$ and $b$ and result $c$, the observed rounding error, $\epsilon$, is given by:
$$c = a \circ b + \epsilon$$
The rounding error of any single operation is extremely small, relative to the architectural precision, but not necessarily small relative to the result of an algebraic addition. Moreover, the rounding errors accumulate and propagate through a sequence of operations. While error propagation depends on the application's dataflow, error introduction depends only on the precision.

We quantify the rounding error for any operation by comparing the observed rounding error (as we describe shortly) to the theoretical maximum units-in-the-last-place ("ulp") error. If the observed error is close to the ulp error, then AMP considers the instruction to be a possible source of numerical faults. We measure closeness to the ulp error using the following Error Ratio:

$$\frac{OBSERVED\ ERROR}{ULP\ ERROR} = \frac{\epsilon}{2^{-(LENGTH\ OF\ MANTISSA)} * c}$$

The Error Ratio is less than or equal to 1 for multiplication/division, but not for addition/subtraction. The ratio is far from 1 if a severe cancellation occurs. We consider cancellation in Section 3.1.3.

Because we would like our profiler to be easily implemented in hardware and the Error Ratio is expensive to compute in hardware, AMP approximates the Error Ratio by scale quantization. Specifically, AMP applies the division operator on just the exponents, resulting in:
$$\frac{2^{\exp(\epsilon)}}{2^{-(LENGTH\ OF\ MANTISSA)} * 2^{\exp(c)}}$$

In log expression, it becomes:

$$\exp(\epsilon) + LENGTH\ OF\ MANTISSA - \exp(c)$$

Threshold: AMP considers a (static) floating point operation to have Large Round-off if a large Error Ratio occurs frequently among dynamic instances of that operation. More precisely, the operation has Large Rounding Error if more than $t_1$% of instances have more than $t_2$% Error Ratio, where $t_1$ and $t_2$ are tunable threshold values.

### 3.1.2 Large Difference in Addend Exponents

Adding two numbers at very different numerical scale (i.e., with large differences in their exponents) will result in the small number being discarded due to the hardware's lack of precision. Consider, the following 4-digit base-10 example: $9.355 \times 10^{30} + 4.267 \times 10^{2}$. With only 4 digits of precision, the result is simply $9.355 \times 10^{30}$.

This behavior is particularly problematic for applications that sum a large number of values. If many of those values are relatively small, every small value may be lost due to this phenomenon, even though the *sum of the small numbers* is large enough to affect the end result [11].

As this behavior occurs regardless of the algebraic signs, it suffices to consider only the absolute values of the numbers. For any addition (or subtraction) with inputs $a$ and $b$, we measure the ratio of the two addends,
$$\frac{|a|}{|b|}, where\ |b| \leq |a|$$
As in Section 3.1.1, the above ratio can be approximated via scale quantization by:
$$|\exp(a) - \exp(b)|$$
Again, we can implement this profile cheaply in hardware by replacing the expensive floating point division with an integer adder and an absolute value unit that are as wide as the exponent.

Threshold: AMP considers a static operation to have a Large Difference in Addend Exponents if more than $t_3$% of instances have exponents that differ by more than $t_4$.

### 3.1.3 Severe Cancellation

Cancellation occurs when two addends are of the same numerical scale but with opposite signs, resulting in the higher order bits being cancelled. Consequently, previously low-order bits get promoted to become the higher order bits of the result. Consider, for example, the base-10 subtraction of $7.65432 \times 10^{4} - 7.6543 \times 10^{4}$. The most significant bit of the result depends on the 6th most significant bit of the inputs. Cancellation is problematic because it can magnify small rounding errors in the previously low-order bits.

Prior work [7][12] has developed checkers to detect cancellation, as well as to quantify the degree of cancellation occurring. Similar to that work, we quantify cancellation by measuring how many of the most significant mantissa bits are cancelled.

Threshold: Because a single cancellation has the potential to ruin an entire application's result, AMP considers an operation to exhibit Severe Cancellation if *any* dynamic instance cancels more than a threshold of mantissa bits, $t_5$.

### 3.1.4 Near Overflow/Underflow

If an operation produces a result that is near the range limit of single precision, then promoting that operation to double precision avoids the possibility of an overflow or underflow.

Threshold: If *any* instance of an operation has a result whose magnitude is greater or smaller than a tunable threshold that is chosen to be near overflow/underflow, then AMP considers that operation to be Near

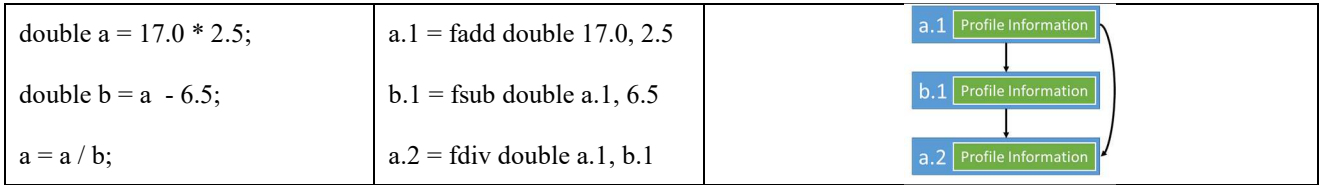

**Figure 1: DDFG annotated with Numerical Behavior**

Overflow/Underflow. A single overflow/underflow is catastrophic and thus AMP cares only whether any instance exceeds this threshold.

## 3.2 Profiler Implementation

Our numerical profiler described thus far works on each floating point operation. We implemented the numerical profiler using a compiler, specifically LLVM (Low Level Virtual Machine) [13]. We developed a compiler pass that instruments the application to emit its numerical profile.

Our compiler pass profiles each operation at the Intermediate Representation (IR) level, instead of at the assembly instruction level. For each IR instruction (henceforth referred to as "instruction") that is a floating point operation, we insert calls to a software library that measure the instruction's numerical behavior. We pass to the software library the input operations to the instruction, as well as a unique identifier for the instruction that consists of the parent function's name and the instruction's destination register. Because LLVM's IR is Single Static Assignment (SSA), each instruction writes to a unique register and we use this name, combined with the name of the function, as an identifier to tie the instruction to its behavior across all instances. We heavily modified SoftFloat[2], a software implementation of an IEEE-754 [14] compliant floating point unit intended for embedded systems, to quantify the numerical behavior of each instruction.

We produce the DDFG annotated with the numerical behaviors as illustrated in Figure 1. The left part of the figure shows the high level code that we are profiling. The center part is the pseudo-LLVM SSA IR to which the code on the left maps. Note that, even though we are writing to the variable "a" twice in the HLL code, the second instance of "a" gets renamed to "a.2" in the IR because of the SSA property. Thus, we label the DDFG on the right with the instruction's name (destination register), thereby allowing for a reverse mapping between the DDFG annotated with the numerical profile and the LLVM IR.

One potential, and subtle, downside of implementing the numerical profiler at the IR level instead of the machine assembly language level, is that the IR may differ from the Instruction Set Architecture (ISA). For example, the IR might contain fused floating point instructions, such as Fused Multiply Add (FMA), not contained in the machine's ISA, or vice versa. This mismatch can cause the IR's DDFG to differ from the assembly language DDFG. We can resolve this issue by ensuring that the Instruction Select phase of the compiler does not fuse or break apart IR instructions. If they are broken apart, the profiling compiler pass has to account for the break. We have not experienced this issue but note that it can occur and is something to be considered when designing the profiling compiler pass.

## 3.3 Profiler Limitations

As with any profile-based static mixed precision scheme, a key limitation is the lack of adaptability to different datasets. For AMP to be useful, the datasets used for profiling or training should be either representative of a class of datasets with which the application is concerned or datasets that are known and used for stress testing numerical behaviors. The mixed precision application produced by AMP, after training on the challenging dataset, should perform well on data that is well-behaved. Choosing and/or developing datasets with bad numerical behaviors is an interesting topic in numerical analysis but beyond the scope of this work.

A limitation that is specific to AMP arises because AMP profiles numerical behavior at the floating point operation level. As such, AMP monitors behavior on a per-operation basis (except for Cancellation) and does not track error propagation in a chain of operations. Some operations might have more of an effect on the final result than others. AMP does not measure the growth of the accumulated error; instead AMP just measures the contributing factors to the error. In this sense, AMP is incomplete.

Despite these limitations, AMP's profiling offers a valuable view into the numerical behavior of an application, especially for complex applications. As we show later, although AMP's profiling is not complete, it suffices for our purposes. Furthermore, AMP has one major advantage with respect to other, more complete schemes discussed in Section 7: cost. By profiling each operation by itself, AMP can be far less costly than schemes that profile the interactions between operations.

## 4. AUTOMATED PRECISION SETTING

AMP takes two steps after profiling: classifying instructions based on their numerical profiles and re-writing the application to promote precision where needed.

## 4.1 Classifying Instructions

Based on the profile it collects, AMP classifies every (IR) instruction in the application into one of four bins. These bins are discussed below in descending priority order; once an instruction has been placed in a bin, it is not considered for any bin further down in the list.

---

[2] http://www.jhauser.us/arithmetic/SoftFloat.html

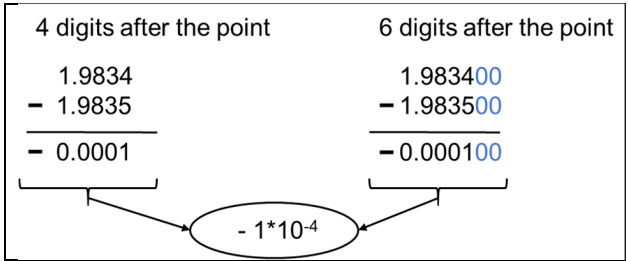

**Figure 2: Greater precision does not reduce cancellation**

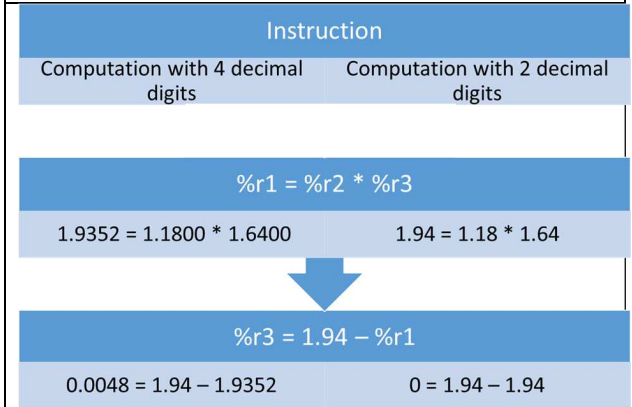

**Figure 3: Promoting instructions in backward slice of instruction with cancellation**

#### 4.1.1 Cancellation Bin

Any instruction that exhibits Severe Cancellation, based on the chosen thresholds, is placed in the Cancellation Bin. Cancellation is different than the other three numerical behaviors profiled by AMP, because of how AMP must set precisions to avoid cancellation. Specifically, promoting the instruction that suffers from cancellation does not help, as we explain later.

#### 4.1.2 Instruction Promotion Bin

Any instruction that exhibits one or more of the other three bad numerical behaviors (Large Round-off, Large Difference in Addend Exponents, or Near Overflow) is placed in the Instruction Promotion Bin. During rewriting, AMP will promote these instructions.

#### 4.1.3 Benign Bin

Instructions in the Benign Bin produce an accurate, rounding-free result every time they are executed with the training dataset. Instructions with no rounding error at the current level of precision will, by definition, not benefit from higher precision.

#### 4.1.4 "Other" Bin

Instructions not placed in any of the bins described above are placed in the "Other" bin. Instructions in this bin may benefit from higher precision but were not explicitly chosen to be in higher precision by the heuristics.

### 4.2 Rewriting the IR

Once AMP has binned all the instructions in the IR, AMP uses two compiler passes—a cancellation cascade pass and a single instruction pass—to promote instructions. Each of these passes promotes certain instructions to higher precision. Because these passes change floating point arithmetic instruction types, we insert cast instructions before and after the promoted instructions as appropriate to ensure that the IR type checker does not fail.

#### 4.2.1 Single Instruction Pass

The single instruction pass applies to all instructions in the Instruction Promotion Bin. AMP promotes each single precision instruction in this bin to be in double precision.

#### 4.2.2 Cancellation Cascade Pass

The cancellation cascade pass operates on instructions in the Cancellation Bin. For an instruction with cancellation, increasing the precision of merely that instruction would not reduce the number of bits being cancelled, as illustrated in Figure 2. The problem is because cancellation is a property of the input operands to the operation and not the operation's precision alone.

Therefore, instead of just increasing the precision of an operation with cancellation, AMP performs backward tracing to reduce the accumulated rounding error in the input operands to that operation. Reducing this rounding error will allow AMP to be more certain that the lower-order bits that get promoted due to the cancellation are more accurate. In Figure 3, we provide an example of this scenario.

For each instruction in the Cancellation Bin, AMP increases its precision and recursively increases the precision of all instructions that produce its operands (i.e., the instruction's backward slice), except those in the Benign Bin and loads. Instructions in the Benign Bin will not benefit from double precision and are a good stopping point for the backwards slice. Our current implementation of AMP does not change the precision of values (data types) in memory nor does it increase the precision of loads.

### 4.3 Choosing Thresholds

The profiler classifies instructions based on the seven threshold values discussed in Section 3.1 and highlighted in Table 1. We consider a chosen set of these values to be a *threshold vector* $\vec{T}$.

Ideally, the programmer would specify the threshold values. However, we currently do not know how to translate these threshold values into "knobs" that a programmer, even one well-versed in numerical analysis,

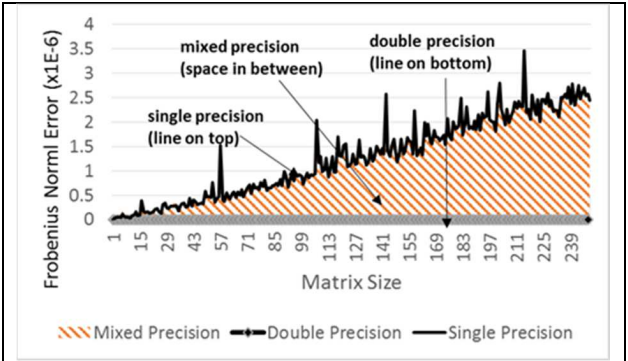

Figure 4: Range of Result Errors for LU with Learning Matrix Size 100x100 drawn from U[-1e6, 1e6]

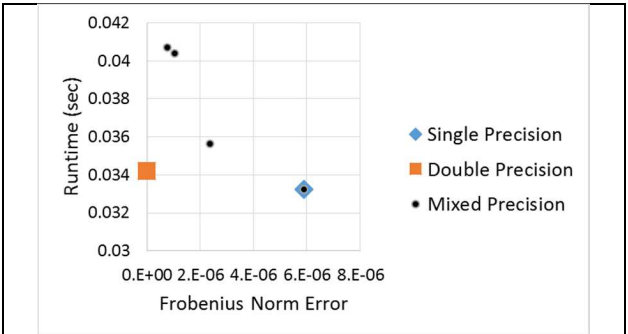

Figure 5: Runtime for Scalar LU. All mixed precision results are strictly worse than double precision.

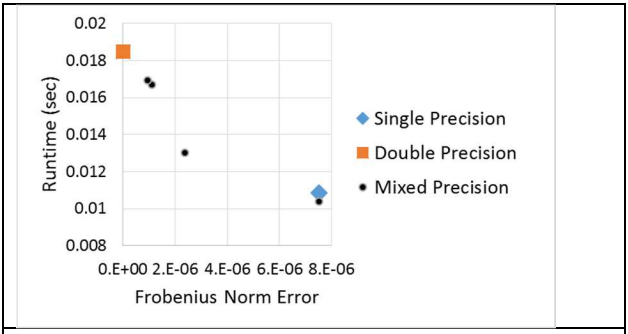

Figure 6: Runtime for Vectorized LU. All mixed precision results offer trade-offs between accuracy and runtime.

would be able to use easily. We leave that mapping or interfacing issue for future work.

In the meanwhile, for each threshold, we sample across a reduced range of quantized values. For each of the 7 thresholds, we look at 6 distinct values in their range leading to $6^7$ sample threshold vectors.

## 5. EXPERIMENTAL EVALUATION

In this section, we use a classic numerical application with well-understood numerical behaviors—LU factorization—as a case study to show the benefits of AMP. Although AMP applies to iterative algorithms, we chose a direct method here in order to decouple the study of numerical behavior with mixed precision from convergence properties. We also evaluate the overhead of AMP's profiling.

### 5.1 Experimental Methodology

We perform all our experiments on a machine running 64-bit Ubuntu with Linux 3.13.0-37-generic. Our machine has an Intel Core i7-3770 processor with 8 MB of cache and 16 GB of main memory. We use LLVM version 3.4 for our automated mixed precision analysis [13]. We vectorize our applications using the automated vectorization pass in LLVM. We repeat all experiments to measure the performance ten times in order to account for variations between the different runs.

### 5.2 Case Study: LU Factorization

The input to this application is a matrix, $A$, that is drawn from different distributions. For a given input matrix, $A$, we factor it into a lower triangular matrix, $L$, and an upper triangular matrix, $U$. For each output matrix, we compare its accuracy against the double precision case using the Frobenius norm. In all experiments, AMP is trained on a 100x100 matrix drawn from the distribution $U[-1e6, 1e6]$.

In Figure 4, we present the accuracy range of the mixed precision results produced by AMP, and we compare these results to single and double precision. The accuracy results are a function of the matrix size. In this graph, for each matrix size, the least error is for double precision (line on bottom), the most error is for single precision (line on top), and mixed precision produces results with errors in the space between double and single precision. For this application, the best mixed precision result produced by AMP does not lead to the double precision result, because LU has large storage arrays whose precision AMP does not modify (as AMP currently only operates on floating point operations). For all metrics, we see four unique mixed precision results for this application given the input dataset.

For performance analysis, we perform experiments on a matrix of size 500x500 after once again training AMP on a matrix of size 100x100. We perform experiments with both scalar code (results in Figure 5) and vectorized code (results in Figure 6). In Figure 5, we see that all of the mixed precision results—with the exception of the one with the same accuracy as single precision—perform worse than double precision. This is because the overheads of the extra cast instructions cannot be hidden in scalar code. In fact, the "best" mixed precision result in terms of accuracy performs significantly worse than double precision.

However, for the vectorized code, every mixed precision result is faster than the double precision result. Thus, for vectorized applications, mixed precision offers

```
// a – a * b
ld.float %r1, [addr1]      // r1 = a
ld.float %r2, [addr2]      // r2 = b
mul.float %r3, %r1, %r2    // r3 = a * b
sub.float %r3, %r1, %r3    // r3 = a – (a * b)
```

|  | $\vec{T_0}$ | $\vec{T_1}$ |
|---|---|---|
| Cancellation Bin | sub.float | N/A |
| Instruction Promotion Bin | N/A | mul.float, sub.float |
| ICS | mul.float, sub.float | mul.float, sub.float |

```
// a – a * b
ld.float %r1, [addr1]         // r1 = a
ld.float %r2, [addr2]         // r2 = b
cast.double %r1, %r1          // r1 = (double) r1
cast.double %r2, %r2          // r2 = (double) r2
mul.double %r3, %r1, %r2      // r3 = a * b
sub.double %r3, %r1, %r3      // r3 = a – (a * b)
```

**Figure 7: Instruction Change Equivalence**

interesting accuracy-versus-performance trade-offs even for small datasets. We do not see a 2X performance difference between double precision and mixed precision due to the cast instructions as well as the bit-shuffle instructions required by mixed precision.

## 5.3 Profiling Costs

Profiling the numerical behavior of applications has two costs. The performance cost is the slowdown due to the calls to the heavily modified SoftFloat library. Recall from Section 3.2 that our software library implements each floating point operation in software, which is the primary cause of the slowdown. The other cost of profiling is the size of the profile itself. The size of the profile depends on the number of floating point operations in the application.

In Table 2, we summarize the costs for LU. We observe, unsurprisingly, that our software implemented profiling incurs very large slowdowns. Fortunately, there are two mitigating factors. First, we have designed our profiler such that it could be easily implemented in hardware, and such hardware would greatly reduce the performance penalty of profiling. Second, we have shown that AMP can—at least in this case study and a few applications not shown here—profile an application on relatively small training datasets (with short running times) and produce a mixed precision version of that application that does well on larger datasets (with long running times).

|  | $\vec{T_0}$ | $\vec{T_1}$ |
|---|---|---|
| Example 1 | <1, 2, 3, 4, 5, **6**> | <1, 2, 3, 4, 5, **7**> |
| Example 2 | <**1**, 2, 3, 4, 5, **6**> | <**2**, 2, 3, 4, 5, **7**> |
| Example 3 | <**1**, 2, 3, 4, 5, **6**> | <**2**, 2, 3, 4, 5, **5**> |

**Figure 8. Prime Threshold Vector Examples. Differences between the two vectors in each example are in bold font.**

The sizes of the profiles are modest. Even profiling LU for over 2 minutes requires a profile that is only 33 MB.

## 6. THRESHOLD VECTOR ANALYSIS

In this section, we delve into the relationships between threshold vectors, in order to shed some light on how vector choices affect outcomes. This analysis is purely for insight; it does not factor into AMP's design or heuristics.

### 6.1 Result Equivalence

One way in which threshold vectors can be related is when multiple threshold vectors lead to the same application result. We call this scenario Result Equivalence. We call a set of threshold vectors that lead to a single result a Result Equivalent Threshold Vector Set, denoted $R$. The size of $R$ can be as small as 1, when there is only one threshold vector that leads to a result. The size of $R$ can be as large as the size of the set of all sample threshold vectors, when there is only one result regardless of the chosen threshold vector, which means the application never benefits from higher precision.

### 6.2 Instruction Change Equivalence

A special case of Result Equivalence is when multiple threshold vectors lead to the precision of the same set of IR instructions being increased. We refer to this scenario as Instruction Change Equivalence (ICE), and we refer to the set of promoted instructions as the Instruction Change Set (ICS). We denote a set of threshold vectors with ICE as $IC$.

We show an example of ICE in Figure 7. In the top third of Figure 7, we present some pseudo-IR code that loads two elements from memory, multiplies them, and then subtracts the product from the first value loaded from memory. In the second third of Figure 7, we have two different threshold vectors in $IC$, $\vec{T_0}$, and $\vec{T_1}$. Using the threshold values in $\vec{T_0}$, AMP places sub.float in the Cancellation Bin. AMP then promotes sub.float and the mul.float that produces one of its inputs. Using the threshold values in $\vec{T_1}$, AMP places sub.float and mul.float in the Instruction Promotion Bin and promotes both of them. Thus, despite $\vec{T_0}$, and $\vec{T_1}$ having different threshold values, AMP ends up promoting the same instructions. The final third of Figure 7 shows the mixed precision version of the pseudo-IR code produced by AMP.

Table 2. Profiling Overheads for LU Decomposition

| Input Dataset Size | Baseline Runtime | Profiled Runtime | Profile Runtime / Baseline Runtime | Profile Size |
|---|---|---|---|---|
| 50 x 50 matrix | 0.759 ms | 15.42 s | 20316 | 4.2 MB |
| 75 x 75 matrix | 0.983 ms | 50.96 s | 51841 | 14 MB |
| 100 x 100 matrix | 1.52 ms | 122.75 s | 80766 | 33 MB |

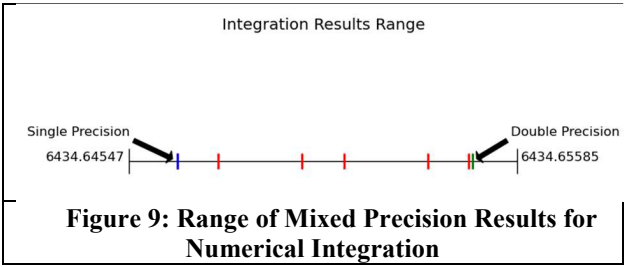

Figure 9: Range of Mixed Precision Results for Numerical Integration

## 6.3 Prime Threshold Vectors

When threshold vectors are related, either by Result Equivalence or Instruction Change Equivalence, we can consider whether some vectors are redundant. That is, we can remove them from $R$ or $IC$ without losing any information or insight.

Consider two threshold vectors in $R$, $\vec{T_0}$, and $\vec{T_1}$, as shown in three examples in Figure 8. In Example 1, $\vec{T_1}$ has the same threshold values as $\vec{T_0}$ except for one threshold value that is larger. Thus, any instruction promotion performed by AMP using $\vec{T_0}$ will also be performed by AMP using $\vec{T_1}$. Therefore, $\vec{T_1}$ is redundant and $\vec{T_0}$ is prime.

In Example 2, $\vec{T_1}$ has the same threshold values as $\vec{T_0}$ except for two threshold values that are larger. Thus, as in Example 1, any instruction promotion performed by AMP using $\vec{T_0}$ will also be performed by AMP using $\vec{T_1}$. Again, $\vec{T_1}$ is redundant and $\vec{T_0}$ is prime.

In Example 3, $\vec{T_1}$ has the same threshold values as $\vec{T_0}$ except for one threshold value that is larger and one that is smaller. Because the vectors differ in different directions on these threshold values, neither vector is redundant and both vectors are prime.

For a given set of vectors with Result Equivalence, we can often gain insight from the set of prime threshold vectors for that set, which we denote as $R'$. When the size of $R'$ is greater than 1, there are relationships between the thresholds, and choosing between two prime vectors is effectively trading off one bad numerical behavior for another. Our results show that, depending on the application and its input dataset, threshold vector equivalence exists and that the size of $R'$ is often greater than one.

## 6.4 Case Study: Numerical Integration

Our case study is numerical integration of $\int_{-10}^{10} \sin(x)\, e^x dx$ using the Gauss-Legendre method [15]. Because this application is not vectorized, we do not expect AMP to provide performance benefits over promoting all variables to double precision; this case study is strictly illustrative and not a good application of AMP.

The accuracy metric that we use for this benchmark is closeness to the double precision result; therefore a result that is closer to the double precision result is said to be more accurate than one that is further.

### 6.4.1 Mixed Precision Results Range

We first show how AMP, depending on the threshold vector, produces different mixed precision configurations with different results. In Figure 9, we present the range of mixed precision results, as well as the single and double precision results. For this benchmark and our set of sample threshold vectors, AMP produces six unique mixed precision results. (Only five are visible in Figure 9, as one of the mixed precision results overlaps with the single precision result.) In this benchmark, none of the mixed precision results equals the double precision result. The reason for this difference is that this benchmark has a small storage array, the precision of which is not increased. The distance between the best mixed precision result and the double precision result is $5 \cdot 10^{-5}$. The difference between the single precision result and the double precision result is $\sim 3.5 \cdot 10^{-3}$. Thus, mixed precision improves the result of this application by two orders of magnitude.

### 6.4.2 Threshold Vector Equivalences

Next, we study how the sampled threshold vectors lead to the set of unique mixed precision results shown in Figure 9. Figure 10 presents the size of $R$ for each of the 6 unique results. The x-axis is the mixed precision results sorted from worst to best. This figure shows that, for our set of sample threshold vectors, the most common result is, with almost a 50% probability, the best result. In Figure 11, we present the size of $R'$ for each of the results. The size of $R'$ can be significantly less than the size of $R$ for all results. For some results, there exist multiple prime threshold vectors while for others there exist just one.

For every prime threshold vector for each result, we look at Instruction Change Equivalence as illustrated in Figure 12. The total height of each bar is the size of $R'$, and the sub-bars are the size of $IC$ for each result. For Result 1, there are five prime threshold vectors with Result Equivalence, leading to four unique ICSs. That is, all four ICSs lead to Result 1. However, Result 2's five prime threshold vectors with Result Equivalence lead to only two unique ICSs. For Results 3, 4, and 5, we have only one prime threshold vector each and thus only one ICS each.

Result 6 has three prime threshold vectors, leading to only one ICS. Although it might seem counter-intuitive that multiple ICSs can lead to the same result, this is possible when we consider that certain instructions in the ICS may not actually affect the result if they are in higher precision. The reason those instructions are allowed to be in higher precision is because we consider each instruction's numerical behavior when making the decision to promote it, not the instruction's impact on the final result.

Figure 13 shows the average percentage of the original IR instructions that are converted to double precision to produce each unique mixed precision result for each threshold vector in $R$. The error bars indicate the range of instructions in double precision for the different ICSs for each result. The x-axis labels are the mixed precision results sorted from worst to best (again with respect to double precision). For this benchmark, the worst mixed precision result (labeled as "Result 1" in Figure 13) has approximately 4% of its IR instructions in double precision with no difference in results compared to keeping all the instructions in single precision. This result means that the precision of this 4% of the IR instructions does not affect the result in any way. The best mixed precision result (labeled as "Result 6"), has around 72% of its IR instructions in double precision. Recall that our mixed precision framework does not achieve the double precision result for this benchmark.

## 7. RELATED WORK

There is related work in several areas, including mixed precision, numerical profiling, and generating numerically challenging datasets.

### 7.1 Programmer-Directed Mixed Precision

A large amount of prior work has explored programmer-directed mixed precision, and this work has been done with both CPUs [16][17][18][19][20][21][22][23], as in our work here, and on GPUs [24][25][26][27][28]. Some other programmer-directed mixed precision differs in that it considers precisions that are lower than the IEEE-754 standard [29][30]. AMP is different from all of this prior work by being automated; it does not require programmer expertise or effort.

### 7.2 Automated Mixed Precision

Like us, some prior researchers have recognized the difficulty of programmer-directed mixed precision and have developed automated tools to perform this task. One class of techniques uses search algorithms to find mixed precision settings [1][2][3]. The automation is helpful, but the drawback is that it performs un-informed searches of the enormous space of mixed precision settings. AMP, in contrast, uses the information from profiling to direct and narrow the search. Another class of techniques uses static analysis of the program to identify where more precision is necessary [31][32][33]. AMP, in contrast has runtime information that can be far more insightful than what can be gleaned from a purely static analysis.

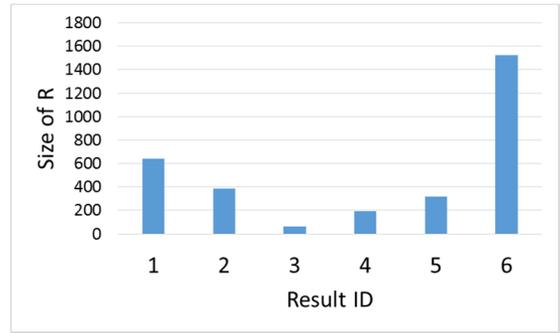

Figure 10: Result Equivalence for Numerical Integration

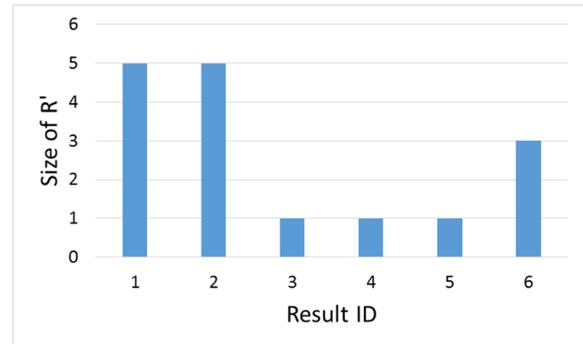

Figure 11: Size of R' for Numerical Integration

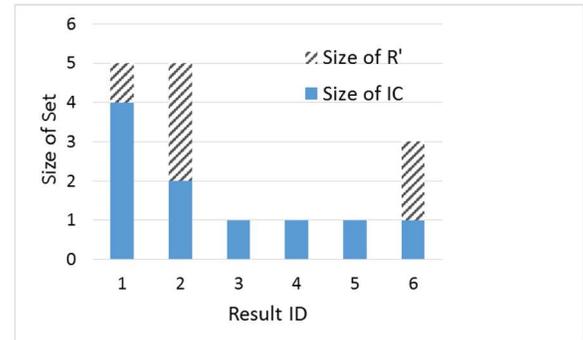

Figure 12: Size of IC vs. R' for Numerical Integration. If size of R' is not visible in bar, then it is same size as IC.

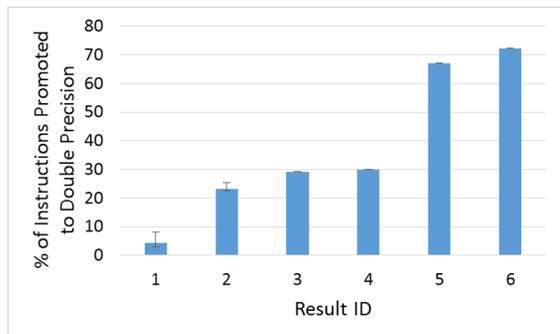

Figure 13: Percentage of Total Instructions in Double Precision for Numerical Integration

## 7.3 Runtime Numerical Profiling

Some prior work [6][7][8][9][10][12][34] has explored runtime profiling of numerical behaviors, which is a key aspect of AMP. These profiling schemes often use binary instrumentation (e.g., with Valgrind) to detect numerical behaviors of interest, although none detect all of the behaviors profiled by AMP. One notable feature of Bao and Zhang [10] is that their scheme focuses on detecting the propagation of a numerical error from one operation to later operations that use this result.

AMP differs from this work in considering a wider range of numerical behaviors and by incorporating numerical profiling into an automated tool that uses the profiles to generate new code. Compared to some prior schemes, AMP does not monitor certain numerical behaviors as precisely (e.g., AMP does not track error across operations), due to our desire for speed and amenability to hardware implementation. However, AMP could incorporate prior profiling schemes if a different trade-off between speed and precision is desired.

## 7.4 Generating Input Datasets

AMP, by virtue of being profile-driven, depends on the training dataset having bad numerical behavior or, at least, worse behavior than the expected input datasets after deployment. Schemes have been developed to generate inputs that produce the most error [35] or that are likely to cause numerically bad behaviors [36]. These schemes are complementary to AMP, in that they could be used to generate training datasets for AMP.

## 8. CONCLUSIONS

In this paper, we present AMP, a profile-driven scheme for automatically setting the precision of each floating point operation. AMP's profiler, which is designed to be easy to implement in hardware, records suspicious numerical behaviors. AMP then decides which operations to promote based on tunable thresholds and automatically rewrites the program at the IR level. Our experimental evaluation shows that AMP can indeed generate a range of mixed precision versions of an application, based on the thresholds, and these different versions can have accuracy/performance tradeoffs. AMP can produce mixed precision applications with better performance than double precision and better accuracy than single precision.

## 9. FUTURE DIRECTIONS

We intend to extend AMP in several ways. First, we would like to extend its applicability. One avenue is to apply AMP to a wider variety of precisions, including half precision and extended precision. Another avenue is applying AMP to iterative numerical methods. An iterative application, such as the Conjugate Gradient method to solve a system of linear equations, runs until it terminates by certain numerical criteria. Mixed precision has been used across iterations for iterative methods at the programmer level (e.g., Baboulin et al. [16]). In this paper, we considered mixed precision computation for a stationary finite-step procedure, which includes a single or a fixed number of iteration steps or a direct (non-iterative) method. Second, to improve confidence in AMP when used on new input datasets, we envision integrating it with prior work that could generate training datasets that lead to particularly bad numerical behavior. Third, we seek develop an interface that allows the programmer to judiciously set the thresholds.

## REFERENCES

[1] C. Rubio-González, C. Nguyen, H. D. Nguyen, J. Demmel, W. Kahan, K. Sen, D. H. Bailey, C. Iancu, and D. Hough, "Precimonious: Tuning Assistant for Floating-point Precision," in *Proceedings of the International Conference on High Performance Computing, Networking, Storage and Analysis*, New York, NY, USA, 2013, pp. 27:1–27:12.

[2] M. O. Lam, J. K. Hollingsworth, B. R. de Supinski, and M. P. Legendre, "Automatically Adapting Programs for Mixed-precision Floating-point Computation," in *Proceedings of the 27th International ACM Conference on International Conference on Supercomputing*, New York, NY, USA, 2013, pp. 369–378.

[3] E. Schkufza, R. Sharma, and A. Aiken, "Stochastic Optimization of Floating-point Programs with Tunable Precision," in *Proceedings of the 35th ACM SIGPLAN Conference on Programming Language Design and Implementation*, New York, NY, USA, 2014, pp. 53–64.

[4] A. Vladimirov, "Arithmetics on Intel's Sandy Bridge and Westmere CPUs: Not All FLOPs Are Created Equal." Apr-2012.

[5] A. Heinecke, K. Vaidyanathan, M. Smelyanskiy, A. Kobotov, R. Dubtsov, G. Henry, A. G. Shet, G. Chrysos, and P. Dubey, "Design and Implementation of the Linpack Benchmark for Single and Multi-node Systems Based on Intel Xeon Phi Coprocessor," in *2013 IEEE 27th International Symposium on Parallel Distributed Processing (IPDPS)*, 2013, pp. 126–137.

[6] A. W. Brown, P. H. J. Kelly, and W. Luk, "Profiling Floating Point Value Ranges for Reconfigurable Implementation," in *Proceedings of the 1st HiPEAC Workshop on Reconfigurable Computing*, 2007, pp. 6–16.

[7] F. Benz, A. Hildebrandt, and S. Hack, "A Dynamic Program Analysis to Find Floating-point Accuracy Problems," in *Proceedings of the 33rd ACM SIGPLAN Conference on Programming Language Design and Implementation*, New York, NY, USA, 2012, pp. 453–462.

[8] D. An, R. Blue, M. Lam, S. Piper, and G. Stoker, "FPInst: Floating Point Error Analysis Using Dyninst."

[9] E. Darulova and V. Kuncak, "Trustworthy Numerical Computation in Scala," in *Proceedings of the 2011 ACM International Conference on Object Oriented Programming Systems Languages and Applications*, New York, NY, USA, 2011, pp. 325–344.

[10] T. Bao and X. Zhang, "On-the-fly Detection of Instability Problems in Floating-point Program Execution," in *Proceedings of the 2013 ACM SIGPLAN International Conference on Object Oriented Programming Systems Languages & Applications*, New York, NY, USA, 2013, pp. 817–832.

[11] T.G. Robertazzi and S. C. Schwartz, "Best 'Ordering' for Floating-Point Addition," *ACM Trans. Math. Softw.*, vol. 14, no. 1, pp. 101–110, Mar. 1988.

[12] M. O. Lam, J. K. Hollingsworth, and G. W. Stewart, "Dynamic Floating-point Cancellation Detection," *Parallel Comput*, vol. 39, no. 3, pp. 146–155, Mar. 2013.

[13] C. Lattner and V. Adve, "LLVM: A Compilation Framework for Life-long Program Analysis and Transformation," in *International Symposium on Code Generation and Optimization*, 2004, pp. 75–86.

[14] "IEEE Standard for Floating-Point Arithmetic," *IEEE Std 754-2008*, pp. 1–58, Aug. 2008.

[15] S. Chapra, *Applied Numerical Methods W/MATLAB: for Engineers & Scientists*, 3 edition. New York: McGraw-Hill Science/Engineering/Math, 2011.

[16] M. Baboulin, A. Buttari, J. Dongarra, J. Kurzak, J. Langou, J. Langou, P. Luszczek, and S. Tomov, "Accelerating Scientific Computations with Mixed Precision Algorithms," *Comput. Phys. Commun.*, vol. 180, no. 12, pp. 2526–2533, 2009.

[17] A. Buttari, J. Dongarra, J. Kurzak, J. Langou, J. Langou, P. Luszczek, and S. Tomov, *Exploiting Mixed Precision Floating Point Hardware in Scientific Computations*. 2007.

[18] A. Buttari, J. Dongarra, J. Kurzak, P. Luszczek, and S. Tomov, *Using Mixed Precision for Sparse Matrix Computations to Enhance the Performance while Achieving 64-bit Accuracy*. 2006.

[19] J. D. Hogg and J. A. Scott, "A Fast and Robust Mixed-precision Solver for the Solution of Sparse Symmetric Linear Systems," *ACM Trans Math Softw*, vol. 37, no. 2, pp. 17:1–17:24, Apr. 2010.

[20] X. S. Li, J. W. Demmel, D. H. Bailey, G. Henry, Y. Hida, J. Iskandar, W. Kahan, S. Y. Kang, A. Kapur, M. C. Martin, B. J. Thompson, T. Tung, and D. J. Yoo, "Design, Implementation and Testing of Extended and Mixed Precision BLAS," *ACM Trans Math Softw*, vol. 28, no. 2, pp. 152–205, Jun. 2002.

[21] R. Strzodka and D. Göddeke, "Mixed precision methods for convergent iterative schemes," *EDGE*, pp. 23–24, 2006.

[22] M. Furuichi, D. A. May, and P. J. Tackley, "Development of a Stokes Flow Solver Robust to Large Viscosity Jumps Using a Schur Complement Approach with Mixed Precision Arithmetic," *J Comput Phys*, vol. 230, no. 24, pp. 8835–8851, Oct. 2011.

[23] D. Göddeke, R. Strzodka, and S. Turek, "Performance and Accuracy of Hardware-oriented Native-, Emulated-and Mixed-precision Solvers in FEM Simulations," *Int J Parallel Emerg Distrib Syst*, vol. 22, no. 4, pp. 221–256, Jan. 2007.

[24] H. Anzt, B. Rocker, and V. Heuveline, "Energy efficiency of mixed precision iterative refinement methods using hybrid hardware platforms," *Comput. Sci. - Res. Dev.*, vol. 25, no. 3–4, pp. 141–148, Aug. 2010.

[25] P. Igounet, E. Dufrechou, M. Pedemonte, and P. Ezzatti, "A Study on Mixed Precision Techniques for a GPU-based SIP Solver," in *2012 Third Workshop on Applications for Multi-Core Architectures (WAMCA)*, 2012, pp. 7–12.

[26] Dominik Goddeke and Robert Strzodka, "Mixed-Precision GPU-Multigrid Solvers with Strong Smoothers," in *Scientific Computing with Multicore and Accelerators*, 0 vols., CRC Press, 2010, pp. 131–147.

[27] S. Le Grand, A. W. Götz, and R. C. Walker, "SPFP: Speed without compromise—A mixed precision model for GPU accelerated molecular dynamics simulations," *Comput. Phys. Commun.*, vol. 184, no. 2, pp. 374–380, Feb. 2013.

[28] D. Göddeke and R. Strzodka, "Cyclic Reduction Tridiagonal Solvers on GPUs Applied to Mixed-Precision Multigrid," *IEEE Trans. Parallel Distrib. Syst.*, vol. 22, no. 1, pp. 22–32, Jan. 2011.

[29] X. Hao and A. Varshney, "Variable-precision Rendering," in *Proceedings of the 2001 Symposium on Interactive 3D Graphics*, New York, NY, USA, 2001, pp. 149–158.

[30] J. Jenkins, E. R. Schendel, S. Lakshminarasimhan, D. A. . I. Boyuka, T. Rogers, S. Ethier, R. Ross, S. Klasky, and N. F. Samatova, "Byte-precision Level of Detail Processing for Variable Precision Analytics," in *Proceedings of the International Conference on High Performance Computing, Networking, Storage and Analysis*, Los Alamitos, CA, USA, 2012, pp. 48:1–48:11.

[31] D. Delmas, E. Goubault, S. Putot, J. Souyris, K. Tekkal, and F. Védrine, "Towards an Industrial Use of FLUCTUAT on Safety-Critical Avionics Software," in *Proceedings of the 14th International Workshop on Formal Methods for Industrial Critical Systems*, Berlin, Heidelberg, 2009, pp. 53–69.

[32] M. D. Linderman, M. Ho, D. L. Dill, T. H. Meng, and G. P. Nolan, "Towards Program Optimization Through Automated Analysis of Numerical Precision," in *Proceedings of the 8th Annual IEEE/ACM International Symposium on Code Generation and Optimization*, New York, NY, USA, 2010, pp. 230–237.

[33] D. Boland and G. A. Constantinides, "A Scalable Precision Analysis Framework," *IEEE Trans. Multimed.*, vol. 15, no. 2, pp. 242–256, Feb. 2013.

[34] E. Tang, E. Barr, X. Li, and Z. Su, "Perturbing Numerical Calculations for Statistical Analysis of Floating-point Program (in)Stability," in *Proceedings of the 19th International Symposium on Software Testing and Analysis*, New York, NY, USA, 2010, pp. 131–142.

[35] W.-F. Chiang, G. Gopalakrishnan, Z. Rakamaric, and A. Solovyev, "Efficient Search for Inputs Causing High Floating-point Errors," in *Proceedings of the 19th ACM SIGPLAN Symposium on Principles and Practice of Parallel Programming*, New York, NY, USA, 2014, pp. 43–52.

[36] E. T. Barr, T. Vo, V. Le, and Z. Su, "Automatic Detection of Floating-point Exceptions," in *Proceedings of the 40th Annual ACM SIGPLAN-SIGACT Symposium on Principles of Programming Languages*, New York, NY, USA, 2013, pp. 549–560.